\documentstyle[amssymb,12pt]{amsart}

\newtheorem{theorem}{Theorem}
\newtheorem{lemma}[theorem]{Lemma}

\newcommand{\forces}{\Vdash}
\newcommand{\res}{\upharpoonright}

\begin{document}

\baselineskip=18pt

  \begin{center}
     {\large Planting Kurepa Trees And Killing Jech--Kunen Trees\\ 
             In a Model By Using One Inaccessible Cardinal}
  \footnote{1980 Mathematics Subject Classification (1985 Revision).
	    Primary 03E35.}
  \end{center}

  \begin{center}
     Saharon Shelah\footnote{The research of the first author
was partially supported by the United States--Israel Binational Science
Foundation, publ. 469.} and Renling Jin
  \end{center}

  \bigskip

  \begin{quote}

    \centerline{Abstract}

    \small

    By an $\omega_{1}$--tree we mean a tree of power $\omega_{1}$
    and height $\omega_{1}\:$. Under $C\!H$ and
    $2^{\omega_{1}}>\omega_{2}$ we call an $\omega_{1}$--tree a
    Jech--Kunen tree if it has $\kappa$ many branches for some $\kappa$
    strictly between $\omega_{1}$ and $2^{\omega_{1}}\:$. 
    In this paper we prove that, assuming the existence of one
    inaccessible cardinal, (1) it is consistent with
    $C\!H$ plus $2^{\omega_{1}}>\omega_{2}$ that there exist
    Kurepa trees and there are no Jech--Kunen trees, which answers
    a question of [Ji2], (2) it is consistent with $C\!H$ plus
    $2^{\omega_{1}}=\omega_{4}$ that only Kurepa trees with $\omega_{3}$
    many branches exist, which answers another question of [Ji2].

  \end{quote}

An partially ordered set, or poset for short, $\langle T,<_{T}\rangle$
is called a tree if for every $t\in T$ the set $\{s\in T:s<_{T}t\}$ is
well--ordered under $<_{T}$. The order type of that set is called the
height of $t$ in $T$, denoted by $ht(t)$. We will not distinguish
a tree from its base set. For every ordinal $\alpha$, let $T_{\alpha}$,
the $\alpha$--th level of $T$, $=\{t\in T:ht(t)=\alpha\}$ and
$T\!\res\!\alpha =\bigcup_{\beta<\alpha}T_{\beta}$. Let $ht(T)$, the height
of $T$, is the smallest ordinal $\alpha$ such that $T_{\alpha}=\emptyset$.
By a branch of $T$ we mean a linearly ordered subset of $T$ which 
intersects every nonempty level of $T$. Let ${\cal B}(T)$ be the set of
all branches of $T$. $T'$ is called a subtree of $T$ if $T'\subseteq T$,
$<_{T'}=<_{T}\bigcap T'\times T'$ ($T'$ inherits the order of $T$) and
for every $\alpha<ht(T')$, $T'_{\alpha}\subseteq T_{\alpha}$.

$T$ is called an $\omega_{1}$--tree if $|T|=\omega_{1}$ and $ht(T)=\omega_{1}$.
An $\omega_{1}$--tree $T$ is called a Kurepa tree if 
$|{\cal B}(T)|>\omega_{1}$ and for every $\alpha\in\omega_{1}$,
$|T_{\alpha}|<\omega_{1}$.
An $\omega_{1}$--tree is called a Jech--Kunen tree if $\omega_{1}<
|{\cal B}(T)|<2^{\omega_{1}}$.

\medskip

T. Jech in [Je1] constructed by forcing a
model of $C\!H$ plus $2^{\omega_{1}}>\omega_{2}$, in which there is a
Jech--Kunen tree. In fact, it is a Kurepa tree with fewer than 
$2^{\omega_{1}}$--many branches. 
Later, K. Kunen [K1] found a model of $C\!H$ plus $2^{\omega_{1}}>\omega_{2}$,
in which there are neither Kurepa trees nor Jech--Kunen trees.
In his paper he gave an equivalent form of Jech--Kunen trees
in terms of compact Hausdorff spaces. The detailed proof can be found in
[Ju, Theorem 4.8].

The second author in [Ji1] started discussing the differences between 
the existence of Kurepa trees and the existence of
Jech--Kunen trees. He showed that it is independent of $C\!H$ plus
$2^{\omega_{1}}>\omega_{2}$ that there exists a Kurepa tree which has
no Jech--Kunen subtrees. He also showed that it is independent of
$C\!H$ plus $2^{\omega_{1}}>\omega_{2}$ that there exists a Jech--Kunen
tree which has no Kurepa subtrees. In his proofs some strongly inaccessible
cardinals were assumed and later, Kunen eliminated the large cardinal
assumption for one of the proofs.

In [SJ], the both authors answered a question of [Ji2] by proving that, 
assuming the existence of
one inaccessible cardinal, it is consistent with $C\!H$ plus 
$2^{\omega_{1}}>\omega_{2}$ that there exist Jech--Kunen trees and
there are no Kurepa trees.

In [Ji2], the second author
proved that, assuming the existence of two inaccessible
cardinals, it is consistent with $C\!H$ plus $2^{\omega_{1}}>\omega_{2}$
that there exist Kurepa trees and there are no Jech--Kunen trees.

Since the consistency of the nonexistence of Jech--Kunen trees implies the
consistency of the existence of an inaccessible cardinal [Ju, Theorem 4.10],
we have to use at least one inaccessible cardinal to build a model of
$C\!H$ plus $2^{\omega_{1}}>\omega_{2}$, in which there are Kurepa trees
but there are no Jech--Kunen trees. The question whether one inaccessible
cardinal is enough was asked in [Ji2]. In this paper, we will give a positive
answer to the question. We also discover that the same techniques can be
used to answer another question in [Ji2] by constructing a model of 
$C\!H$ plus $2^{\omega_{1}}=\omega_{4}$, in which only the Kurepa trees
with $\omega_{3}$--many branches exist.

First let's look at the second author's original idea in [Ji2] to construct a
model of $C\!H$ plus $2^{\omega_{1}}>\omega_{2}$, in which there are
Kurepa trees and there are no Jech--Kunen trees, by using two inaccessible
cardinals. Let $\kappa_{1} <\kappa_{2}$ be two strongly inaccessible
cardinals in a model $M$. First, Jin collapses $\kappa_{2}$ down to
$\kappa_{1}^{+}$ by forcing with a $<\kappa_{1}$--support
L\'{e}vy collapsing order. 
Next, he collapses $\kappa_{1}$ down to
$\omega_{1}$ by forcing with a finite support L\'{e}vy collapsing order.
This step will create a very homogeneous Kurepa tree. Then he force with
that Kurepa tree $\lambda$--many times for some regular cardinal $\lambda$
which is greater than $\kappa_{2}$. In the resulting model, that Kurepa tree
has $\lambda$--many branches and $\lambda =2^{\omega_{1}}$. In that model
there are no Jech--Kunen trees.

If we want to obtain the same result by using only one inaccessible cardinal,
we may have to find a way to create a homogeneous 
$\omega_{1}$--tree with every level 
countable, without the assistance of large cardinals. Then the questions
arise. First, how can we create the desired tree? Second, can we force with
that tree for multiple times (with countable supports) without collapsing 
$\omega_{1}$. (Note that a normal $\omega_{1}$--tree with every level
countable is never $\omega_{1}$--closed.)

In this paper, we construct a homogeneous generic $\omega_{1}$--tree with every 
level countable
by forcing with an $\omega_{1}$--closed poset, whose elements are countable
homogeneous normal subtrees of $\langle 2^{<\omega_{1}},\subseteq\rangle$. 
The generic tree is, in fact, a Suslin tree.
Then we force with that generic tree $\lambda$--many times with countable
supports. We will prove that this two--step forcing adds no new countable
sequences of ordinals, hence it will not collapse $\omega_{1}$. We will 
also prove that if the ground model is Silver's model (see [K2, pp. 259]),
then in the final model there are no Jech--Kunen trees.

Before proving our results we need more notations and definitions.

\medskip
 
A tree $T$ is called normal if,

(1) every $t\in T$, which is not in the top level of $T$, 
has at least two immediate successors,

(2) for every limit ordinal $\alpha<ht(T)$ and every 
$B\in {\cal B}(T\!\res\!\alpha)$, there exist at most one least upper
bound of $B$ in $T$,

(3) for every $t\in T$ and $\alpha$ such that $ht(t)<\alpha <ht(T)$, there
exists $t'\in T_{\alpha}$ such that $t<_{T} t'$.

Without loss of generality, we consider only the trees which are subtrees of 
$\langle 2^{<\omega_{1}},\subseteq\rangle$ with the unique root $\emptyset$.
Let $T$ be a tree and $B\subseteq T$ be a totally ordered subset of $T$.
Then $\bigcup B$ is the only candidate for the least upper bound of $B$ in $T$.

Let $\alpha\in\omega_{1}$ and $s,t\in 2^{\alpha}$. We define a map $F_{s,t}$
from $2^{<\omega_{1}}$ to $2^{<\omega_{1}}$. Let $u\in 2^{\beta}$ for some
$\beta <\omega_{1}$. The domain of $F_{s,t}(u)$ is $\beta$ and
for every $\gamma<\beta$, if $\gamma<\alpha$, then let\[F_{s,t}(u)(\gamma)
=u(\gamma)+t(\gamma)-s(\gamma)\;\;(\mbox{mod 2}),\] otherwise let $F_{s,t}(u)
(\gamma)=u(\gamma)$.

\begin{lemma}

$F_{s,t}(s)=t$, $F_{s,t}(t)=s$ and $F_{s,t}$ is an isomorphism from
$\langle 2^{<\beta},\subseteq\rangle$ to $\langle 2^{<\beta},\subseteq\rangle$
for any $\beta\leq\omega_{1}$.

\end{lemma}

\noindent {\bf Proof:}\quad Trivial. \quad $\Box$

\bigskip

A normal tree $T$ is called homogeneous if for any $\alpha <ht(T)$, for any
$s,t\in T_{\alpha}$, $F_{s,t}\!\res\! T$ is an isomorphism from $T$ to $T$.

Let \[{\Bbb P}_{hom}=\{T:T\mbox{ is a countable homogeneous normal subtree
of }\langle 2^{<\omega_{1}},\subseteq\rangle\}\]
be a poset ordered by letting $T<T'$ iff $ht(T')<ht(T)$ and $T'=T\!\res\!
ht(T')$.

\begin{lemma}

Let ${\cal T}$ be a totally ordered subset of ${\Bbb P}_{hom}$. Then
$\bigcup {\cal T}$ is a homogeneous tree. Moreover, if ${\cal T}$ is countable,
then $\bigcup {\cal T}\in {\Bbb P}_{hom}$.

\end{lemma}

\noindent {\bf Proof:}\quad Trivial. \quad $\Box$

\bigskip

\noindent {\bf Remark:}\quad Above lemma says that ${\Bbb P}_{hom}$ is 
$\omega_{1}$--closed, which means that every countable decreasing
sequence in ${\Bbb P}_{hom}$ has a lower bound in ${\Bbb P}_{hom}$.

\begin{lemma}

Let $T\in {\Bbb P}_{hom}$ and $ht(T)=\alpha$ for some limit ordinal
$\alpha<\omega_{1}$. Let ${\cal C}$ be a countable subset of 
${\cal B}(T)$. Then there exists $\overline{T}\in {\Bbb P}_{hom}$
such that $\overline{T}<T$ and for every $C\in {\cal C}$,
$\bigcup C\in\overline{T}_{\alpha}$.

\end{lemma}

\noindent {\bf Proof:}\quad Without loss of generality, we can assume 
that for every $t\in T$, there exists $C\in {\cal C}$ such that
$t\in C$. (This will guarantee the normality of the resulting tree.)
We now construct inductively a sequence of countable trees
$\langle T_{n}:n\in\omega\rangle$ such that:

(1) $T_{0}=T\bigcup\{\bigcup C: C\in {\cal C}\}$,

(2) for every $n\in\omega$, $ht(T_{n})=\alpha +1$ and

(3) for every $n\in\omega$, \[ T_{n+1}=T_{n}\bigcup\{F_{s,t}(u):
s,t\in T_{n}, \; ht(s)=ht(t)\mbox{ and }u\in (T_{n})_{\alpha}\}.\]

Note that if $I$ is an isomorphism from $T$ to $T$,
then for every $B\in {\cal B}(T)$, $I[B]\in {\cal B}(T)$.

Let $\overline{T}=\bigcup_{n\in\omega}T_{n}$. It is obvious that
$\overline{T}$ is countable and for any $s,t\in \overline{T}$ such that
$ht(s)=ht(t)$, $F_{s,t}$ is an isomorphism from $\overline{T}$ to
$\overline{T}$. Hence $\overline{T}\in {\Bbb P}_{hom}$, $\overline{T}
<T$ and for every $C\in {\cal C}$, $\bigcup C\in T_{0}\subseteq
\overline{T}$. \quad $\Box$

\bigskip

Next we discuss forcing method.
For the terminology and basic facts of forcing, see [K2] and [Je2].
We always assume the consistency of $Z\!F\!C$ and let $M$ be always a
countable transitive model of $Z\!F\!C$. In the forcing arguments,
we always let $\dot{a}$ be a name of $a$. 
For every element $a$ in the ground model, we may use $a$ itself as its name.
Let $\Bbb P$ be a poset in a model $M$, $\dot{a}$ be a $\Bbb P$--name for $a$
and $G$ be a $\Bbb P$--generic filter over $M$. Then $\dot{a}_{G}$ 
is the value of $\dot{a}$ in $M[G]$ (see [K2, pp. 189] for the definition of
$\dot{a}_{G}$). 

Let $I,J$ be two sets. Let 
\[Fn(I,J,\omega_{1})=\{p:p\subseteq I\times J
\mbox{ is a function and }|p|<\omega_{1}\}\] 
be a poset ordered by
reverse inclusion. Let $I$ be a subset of a cardinal $\kappa$. Let
\[Lv(I,\omega_{1})=\]
\[\{p:p\subseteq (I\times\omega_{1})\times\kappa
\mbox{ is a function, }|p|<\omega_{1}\mbox{ and }
\forall\langle\alpha,\beta\rangle\in\mbox{dom}(p)
(p(\alpha,\beta)\in\alpha)\}\]
be a poset ordered by reverse inclusion. The poset $Lv(\kappa,\mu)$
for some regular cardinals $\kappa>\mu$ is usually called
a $<\mu$--support L\'{e}vy collapsing order.
Let $T$ be a tree and $I$ be an index set. For a function $p$ from $I$ to $T$,
let $supt(p)$, the support of $p$, be the set $\{i\in I:p(i)\neq\emptyset\}$. 
Let \[{\Bbb P}(T,I,\omega_{1})
=\{p:p\in T^{I},\;|supt(F)|<\omega_{1}\}.\]
For any $p,p'\in {\Bbb P}
(T,I,\omega_{1})$, define $p\leq p'$ iff for every $i\in I$, 
$p'(i)\leq_{T} p(i)$. Let $\Bbb R$ be a poset and $\dot{T}$ be an 
$\Bbb R$--name for a tree $T$. Let
\[{\Bbb P}(\dot{T},I,\omega_{1})=\{\dot{q}:\dot{q}\in (\dot{T})^{I},\;
|supt(\dot{q})|<\omega_{1}\}.\]
Then ${\Bbb P}(\dot{T},I,\omega_{1})$ is an $\Bbb R$--name for the poset
${\Bbb P}(T,I,\omega_{1})$.
Let ${\Bbb Q}={\Bbb P}(T,I,\omega_{1})$ (or ${\Bbb P}(\dot{T},I,\omega_{1})$)
and $J\subseteq I$. We denote ${\Bbb Q}\!\res\! J$ for the set
$\{p\!\res\! J:p\in {\Bbb Q}\}$. If $H$ is a $\Bbb Q$--generic filter, then
let $H_{J}=\{p\!\res\! J:p\in H\}$.

\begin{lemma}

Let $T$ be an $\omega_{1}$--tree and ${\Bbb P}$ be an $\omega_{1}$--closed
poset in a model $M$. Let $G$ be a $\Bbb P$--generic filter over $M$.
Assume that there exists a branch of $T$ in $M[G]\smallsetminus M$. Then
$T$ is neither a Kurepa tree nor a Jech--Kunen tree in $M$. Moreover,
there exists an isomorphic embedding from $\langle 2^{<\omega_{1}},\subseteq
\rangle$ into $T$.

\end{lemma}

\noindent {\bf Proof:}\quad See [K2, pp. 260] and [Ju, Theorem 4.9].\quad $\Box$

\begin{lemma}

Let $M$ be a model, ${\Bbb P}=({\Bbb P}_{hom})^{M}$ and $G$ be a 
$\Bbb P$--generic filter over $M$. Let $T_{G}=\bigcup G$. Then
the generic tree $T_{G}$ is a homogeneous normal $\omega_{1}$--tree
with every level countable. In fact, $T_{G}$ is a Suslin tree.

\end{lemma}

\noindent {\bf Proof:}\quad See [Je2, Theorem 48] for the proof.
The homogeneity of $T_{G}$ follows from Lemma 2. \quad $\Box$

\begin{lemma}

Let $M$ be a model, $I$ be an index set in $M$, ${\Bbb P}=
({\Bbb P}_{hom})^{M}$, $T_{\dot{G}}$ be $\Bbb P$--name for the
$\Bbb P$--generic tree $T_{G}$, and $\dot{\Bbb Q}=
{\Bbb P}(T_{\dot{G}},I,\omega_{1})$,
which is a $\Bbb P$--name for ${\Bbb P}(T_{G},I,\omega_{1})$. Let
$G*H$ be a ${\Bbb P}*\dot{\Bbb Q}$--generic filter over $M$. Then
$M^{\omega}\bigcap M[G*H]\subseteq M$.

\end{lemma}

\noindent {\bf Proof:}\quad
Suppose that there is an $f\in M^{\omega}\bigcap M[G*H]$ such that
$f\not\in M$. Let $\langle p,\dot{q}\rangle\in {\Bbb P}*\dot{\Bbb Q}$
such that \[\langle p,\dot{q}\rangle\forces \dot{f}\in A^{\omega}
\smallsetminus M\] for some $A\in M$.

We now want to construct a sequence $\langle\langle p_{n},
\dot{q}_{n}\rangle\in {\Bbb P}*\dot{\Bbb Q}: n\in\omega\rangle$ in $M$
such that for every $n\in \omega$,

(1) $\langle p_{n+1},\dot{q}_{n+1}\rangle\leq
\langle p_{n},\dot{q}_{n}\rangle\leq\langle p,\dot{q}\rangle$,

(2) $\exists a_{n}\in A\;(\langle p_{n},\dot{q}_{n}\rangle
\forces\dot{f}(n)=a_{n})$,

(3) $\forall i\in supt(\dot{q}_{n})\;\exists t_{n}(i)\in
p_{n}\smallsetminus p_{n-1}\;(p_{n}\forces\dot{q}_{n}(i)
=t_{n}(i))$.

The contradiction follows from the construction. Let $\overline{p}_
{\omega}=\bigcup_{n\in\omega}p_{n}$. For every
$i\in\bigcup_{n\in\omega}supt(\dot{q}_{n})$, let \[C_{i}=\{t\in\overline{p}_
{\omega}:\exists n\in\omega,\;t<t_{n}(i)\}.\] By (3),
$C_{i}\in {\cal B}(\overline{p}_{\omega})$. By Lemma 3, there is $p_{\omega}
\in {\Bbb P},\; p_{\omega}
\leq\overline{p}_{\omega}$ such that $\bigcup C_{i}\in p_{\omega}$.
Define $\dot{q}_{\omega}$ from $I$ to $T_{\dot{G}}$ such that
$\dot{q}_{\omega}(i)=\bigcup C_{i}$ if $i\in\bigcup_{n\in\omega}
supt(\dot{q}_{n})$ and $\dot{q}_{\omega}(i)=\emptyset$ otherwise.
(In fact, $q$ is in $M$.)
Then $\langle p_{\omega},\dot{q}_{\omega}\rangle\leq \langle p_{n},
\dot{q}_{n}\rangle$ for every $n\in\omega$. Hence, for every $n\in\omega$,
\[\langle p_{\omega},\dot{q}_{\omega}\rangle\forces\dot{f}(n)=a_{n}.\]
This contradicts $f\not\in M$.

Assume that we have already had $\langle p_{n},\dot{q}_{n}\rangle$
for every $n\leq m$.

First, let $\langle r,\dot{s}\rangle\leq\langle p_{m},\dot{q}_{m}\rangle$
and $a_{m+1}\in A$ such that
\[\langle r,\dot{s}\rangle\forces\dot{f}(m+1)=a_{m+1}.\]
For every $i\in supt(\dot{s})$, \[r\forces\exists\alpha\in\omega_{1}\;
(\dot{s}(i)\in 2^{\alpha}).\]
Then there exist $\alpha\in\omega_{1}$ and $r'\leq r$ such that
\[r'\forces \dot{s}(i)\in 2^{\alpha}.\]
Since $\Bbb P$ is $\omega_{1}$--closed and 
\[r'\forces\mbox{ The domain of }\dot{s}(i)\mbox{ is }\alpha,\mbox{ a countable
ordinal.}\] 
then there exist $t(i)\in 2^{\alpha}$ and $r''\leq r'$ such that
\[ r''\forces\dot{s}(i)=t(i).\]
Let $r'''\leq r''$ such that $ht(r''')>\max\{\alpha,ht(p_{m})\}$. Then
\[r'''\forces\dot{s}(i)=t(i)\in r'''\]
because $\forces\dot{s}(i)\in T_{\dot{G}}$.

Since $supt(\dot{s})$ is countable and $\Bbb P$ is $\omega_{1}$--closed,
then we can find $p_{m+1}\leq r'''$ such that
\[\forall i\in supt(\dot{s})\;\exists\alpha<ht(p_{m+1})\;\exists t(i)\in
(p_{m+1})_{\alpha}\;(p_{m+1}\forces\dot{s}(i)=t(i)).\]
Let $t_{m+1}(i)\in p_{m+1}\smallsetminus p_{m}$ such that $t(i)\leq
t_{m+1}(i)$ and define $\dot{q}_{m+1}(i)=t_{m+1}(i)$ if
$i\in supt(\dot{s})$ and $\dot{q}_{m+1}(i)=\emptyset$ otherwise.
This ends the construction and the sequence we have constructed does
obviously satisfy (1), (2) and (3). \quad $\Box$

\bigskip

\noindent {\bf Remark:}\quad The poset ${\Bbb P}*\dot{\Bbb Q}$ in Lemma 6
is, in fact, strategically complete. Let $\Bbb R$ be any poset.
Two players, $I$ and $I\!I$, choose from $\Bbb R$ successively the members of
a decreasing sequence
\[a_{0}\geq b_{0}\geq a_{1}\geq b_{1}\geq\cdots
\geq a_{n}\geq b_{n}\geq\cdots.\] 
$I$ chooses the $a_{n}$'s and $I\!I$ chooses the $b_{n}$'s. $I\!I$ wins the
game if and only if the sequence has a lower bound in $\Bbb R$. We call 
$\Bbb R$ strategically complete if $I\!I$ has a winning strategy.
It has been shown that $\Bbb R$ is strategically complete if and only if
there exists a poset $\Bbb S$ such that ${\Bbb R}\times {\Bbb S}$ has a dense
subset which is $\omega_{1}$--closed (see [Je3, pp. 90]).

\begin{theorem}

Assuming the existence of an inaccessible cardinal, it is consistent with
$C\!H$ plus $2^{\omega_{1}}>\omega_{2}$ that there exist Kurepa tree
but there are no Jech--Kunen trees.

\end{theorem}

\noindent {\bf Proof:}\quad
Let $M$ be a model of $G\!C\!H$, $\kappa$ be an inaccessible cardinal
and $\lambda>\kappa$ be a regular cardinal in $M$.
In $M$, let ${\Bbb P}_{1}=Lv(\kappa,\omega_{1})$, ${\Bbb P}_{2}=
{\Bbb P}_{hom}$, $T_{\dot{G_{2}}}$ be a ${\Bbb P}_{2}$--name for the
${\Bbb P}_{2}$--generic tree $T_{G_{2}}=\bigcup G_{2}$, where
$G_{2}$ is a ${\Bbb P}_{2}$--generic filter, and $\dot{\Bbb Q}=
{\Bbb P}(T_{\dot{G}_{2}},\lambda,\omega_{1})$.
Let $G_{1}\times (G_{2}*H)$ be a ${\Bbb P}_{1}\times 
({\Bbb P}_{2}*\dot{\Bbb Q})$--generic filter over $M$.
We will show that $M[G_{1}\times (G_{2}*H)]=M[G_{1}][G_{2}*H]$ 
is the model we are looking for.

\bigskip

{\bf Claim 7.1}\quad $M^{\omega}\bigcap M[G_{1}][G_{2}*H]
\subseteq M$.

{\bf Proof of Claim 7.1 :}\quad
By Lemma 6, $M^{\omega}\bigcap M[G_{2}*H]\subseteq M$. This implies that
${\Bbb P}_{1}$ is still $\omega_{1}$--closed in $M[G_{2}*H]$.
Hence $(M[G_{2}*H])^{\omega}\bigcap M[G_{2}*H][G_{1}]\subseteq
M[G_{2}*H]$. So for every $f\in M^{\omega}\bigcap M[G_{2}*H][G_{1}]$,
$f$ is in $M[G_{2}*H]$ and hence, $f$ is in $M$. The claim is true
because $M[G_{1}][G_{2}*H]=M[G_{2}*H][G_{1}]$.

\bigskip

{\bf Claim 7.2}\quad ${\Bbb P}_{1}\times ({\Bbb P}_{2}*\dot{\Bbb Q})$
has the $\kappa$--c.c..

{\bf Proof of Claim 7.2 :}\quad
A poset $\Bbb R$ is called $\lambda$--centered for some regular cardinal
$\lambda$ if for any subset $S\subseteq 
{\Bbb R}$ and $|S|\geq\lambda$, there exists $S'\subseteq S$, 
$|S'|\geq\lambda$, such that any two elements in $S'$ are compatible. 
By a simple $\Delta$--system lemma argument, we can show that ${\Bbb P}_{1}$
is $\kappa$--centered. Since $|{\Bbb P}_{2}|=\omega_{1}$, then
$|T_{\dot{G}_{2}}|\leq (|{\Bbb P}_{2}|^{\omega_{1}})^{\omega_{1}}=\omega_{2}$.
Again by a simple $\Delta$--system lemma argument, we can show that
${\Bbb P}_{2}*\dot{Q}$ is $\kappa$--centered. In fact, it is 
also $\omega_{3}$--centered. Hence
${\Bbb P}_{1}\times ({\Bbb P}_{2}*\dot{Q})$ is $\kappa$--centered, which
implies the $\kappa$--c.c..

\bigskip

\noindent {\bf Remark:}\quad By Claim 1 and Claim 2 and the fact that
$M[G_{1}]\models [C\!H + 2^{\omega_{1}}=\omega_{2}=\kappa]$, 
we know that $\omega_{1}$ and all the cardinals greater than or equal to
$\kappa$ in $M$ is preserved in $M[G_{1}][G_{2}*H]$. We also know that
$M[G_{1}][G_{2}*H]\models [C\!H + 2^{\omega_{1}}=\lambda >\kappa]$.

\bigskip

{\bf Claim 7.3}\quad
$T_{G_{2}}$ is a Kurepa tree with $\lambda$--many branches in
$M[G_{1}][G_{2}*H]$.

{\bf Proof of Claim 7.3 :}\quad
It is obvious that $T_{G_{2}}$ is an $\omega_{1}$--tree with every level
countable (in fact, it is a Suslin tree in $M[G_{2}]$). In $M[G_{1}][G_{2}]$,
$\dot{Q}_{G_{2}}={\Bbb P}(T_{G_{2}},\lambda,\omega_{1})$ is a countable
support (note that no new countable sequences of ordinals are added) product of
$\lambda$--many copies of $T_{G_{2}}$. Then forcing with $\Bbb Q$ will add
at least $\lambda$--many new branches to $T_{G_{2}}$. Hence
$\lambda\leq |{\cal B}(T_{G_{2}})|\leq 2^{\omega_{1}}=\lambda$.

\bigskip

{\bf Claim 7.4}\quad
There are no Jech--Kunen trees in $M[G_{1}][G_{2}*H]$.

{\bf Proof of Claim 7.4 :}\quad
Suppose that $T$ is a Jech--Kunen tree in $M[G_{1}][G_{2}*H]$.

Since $|T|=\omega_{1}$, then there exists a cardinal $\theta<\kappa$
and a subset $I$ of $\lambda$ with $|I|\leq\omega_{2}$
such that $T\in M[G'_{1}][G_{2}*H_{I}]$, where $G'_{1}=G_{1}\bigcap
Lv(\theta,\omega_{1})$ and $H_{I}=H\bigcap {\Bbb Q}\!\res\! I$. 
This is true because 
${\Bbb P}_{1}$ has the $\kappa$--c.c. and ${\Bbb P}_{2}*\dot{\Bbb Q}$
has the $\omega_{3}$--c.c.. 
In $M[G'_{1}][G_{2}*H_{I}]$, $2^{\omega_{1}}<\kappa$, so that there exists
a branch $b$ of $T$ in $M[G_{1}][G_{2}*H]\smallsetminus 
M[G'_{1}][G_{2}*H_{I}]$. Since $Lv(\kappa\smallsetminus\theta,\omega_{1})$
in $M$ is still $\omega_{1}$--closed in $M[G'_{1}][G_{2}*H_{I}]$ 
and $T$ is a Jech--Kunen tree in 
$M[G_{1}][G_{2}*H]$, then by Lemma 4, $b\not\in M[G_{1}][G_{2}*H_{I}]$.

Let $M'=M[G_{1}][G_{2}*H_{I}]$. We now work in $M'$. 
In $M'$, ${\Bbb Q}\!\res\!(\lambda
\smallsetminus I)$ has the $\omega_{1}$--c.c.. Then there exists $J\subseteq
\lambda\smallsetminus I$ with $|J|=\omega_{1}$ in $M'$ such that
$b\in M'[H_{J}]$.
Let $r\in H_{J}$ be such that
\[r\forces_{{\Bbb Q}\!\res\! J}\exists b\in {\cal B}(T)\smallsetminus M'.\]
Since $T_{G_{2}}$ is homogeneous (here we use the homogeneity of the tree), 
then we can assume that \[\forces_{{\Bbb Q}\!\res\! J}\exists b\in
{\cal B}(T)\smallsetminus M'.\]
By the maximal principle we can find a ${\Bbb Q}\!\res\! J$--name $\dot{b}$
for $b$ in $M'$ such that \[\forces_{{\Bbb Q}\!\res\! J}\dot{b}\in {\cal B}(T)
\smallsetminus M'.\]
Since $b\not\in M'$, then in $M'$, the sentence $\Phi ({\Bbb Q}\!\res\! J,
T,\dot{b})$ is true, where $\Phi (X,Y,Z)$ is \[\forall s\in X\;
\exists s_{0},s_{1}\leq s\;\exists\alpha\in\omega_{1}\;\exists t_{0},t_{1}
\in Y_{\alpha},\; t_{0}\neq t_{1},\;(s_{i}\forces t_{i}\in Z)\mbox{ for }
i=0,1.\]
In $M[G_{1}][G_{2}*H]$ $T$ has less than $\lambda$--many branches, so
there exists $\mu<\lambda$ such that $I\bigcup J\subseteq\mu$ and
every branch of $T$ in $M[G_{1}][G_{2}*H]$ is already in
$M'[H_{\mu\smallsetminus I}]$.
Let $J'\subseteq\lambda\smallsetminus\mu$ be such that $|J'|=|J|$ and
let $\pi$ be the natural isomorphism from ${\Bbb Q}\!\res\! J$ to
${\Bbb Q}\!\res\! J'$. Then in $M'$
\[\forces_{{\Bbb Q}\!\res\! J'}\pi_{*}(\dot{b})\in {\cal B}(T)\]
is true and \[M'\models\Phi ({\Bbb Q}\!\res\! J',T,\pi_{*}(\dot{b})),\]
where $\pi_{*}$ is the map from ${\Bbb Q}\!\res\! J$--names to
${\Bbb Q}\!\res\! J'$--names induced by $\pi$ (see [K2, pp. 222] for
the definition of $\pi_{*}$). 

\medskip

{\bf Subclaim 7.4.1}\quad $M'[H_{\mu\smallsetminus I}]\models
[\Phi ({\Bbb Q}\!\res\! J',T,\pi_{*}(\dot{b}))$ and $\forces_
{{\Bbb Q}\!\res\! J'}\pi_{*}(\dot{b})\in {\cal B}(T)]$.

{\bf Proof of Subclaim 7.4.1 :}\quad
Let $H_{J'}$ be a ${\Bbb Q}\!\res\! J'$--generic filter over 
$M'[H_{\mu\smallsetminus I}]$. Then $H_{J'}$ is also a 
${\Bbb Q}\!\res\! J'$--generic filter over $M'$.
Hence in $M'[H_{J'}]$, $(\pi_{*}(\dot{b}))_{H_{J'}}\in {\cal B}(T)$.
If $s_{i}\in H_{J'}$, then $t_{i}\in (\pi_{*}(\dot{b}))_{H_{J'}}$ is also
true in $M'[H_{J'}]$. 

In $M'[H_{J'}]$, forcing with ${\Bbb Q}\!\res\!(\mu\smallsetminus I)$ will not
change the truth of the above sentences. Hence in $M'[H_{J'}][H_
{\mu\smallsetminus I}]=M'[H_{\mu\smallsetminus I}][H_{J'}]$, 
$(\pi_{*}(\dot{b}))_{H_{J'}}\in {\cal B}(T)$ and
$t_{i}\in (\pi_{*}(\dot{b}))_{H_{J'}}$ are also true.
This implies that \[ M'[H_{\mu\smallsetminus I}]\models
[\Phi ({\Bbb Q}\!\res\! J',T,\pi_{*}(\dot{b}))\mbox{ and }
\forces_{{\Bbb Q}\!\res\! J'}\pi_{*}(\dot{b})\in {\cal B}(T)].\]
This ends the proof of Subclaim 7.4.1.

\medskip
 
Since forcing with ${\Bbb Q}\!\res\! J'$
will not add any new branches of $T$, then $B=(\pi_{*}(\dot{b}))_{H_{J'}}$ 
is already in $M'[H_{\mu\smallsetminus I}]$. In $M'[H_{\mu\smallsetminus I}]$,
let \[D=\{r\in {\Bbb Q}\!\res\! J':
\exists t\not\in B\;(r\forces_{{\Bbb Q}\!\res\! J'}
t\in\pi_{*}(\dot{b}))\}.\] Then
$D$ is dense in ${\Bbb Q}\!\res\! J'$ because $\Phi ({\Bbb Q}\!\res\! J',
T,\pi_{*}(\dot{b}))$ is true in $M'[H_{\mu\smallsetminus I}]$. 
If $r_{0}\in D\bigcap H_{J'}$, then
$r_{0}\forces\pi_{*}(\dot{b})\neq B$. This contradicts
$(\pi_{*}(\dot{b}))_{H_{J'}}=B$. \quad $\Box$

\begin{theorem}

Assuming the existence of an inaccessible cardinal, it is consistent with
$C\!H$ plus $2^{\omega_{1}}=\omega_{4}$ that only the Kurepa trees with
$\omega_{3}$--many branches exist. 

\end{theorem}
 
\noindent {\bf Proof:}\quad
Let's follow the notation of the proof of Theorem 7.
Let $\lambda=\kappa^{+}$ in $M$. Let \[{\Bbb P}_{3}= 
Fn(\kappa^{++},2,\omega_{1})=Fn(\omega_{4},2,\omega_{1})\] in 
$M[G_{1}][G_{2}*H]$ (note that ${\Bbb P}_{3}$ is absolute with respect to
$M$ and $M[G_{1}][G_{2}*H]$). Let $G_{3}$ be a ${\Bbb P}_{3}$--generic filter
over $M[G_{1}][G_{2}*H]$. In $M[G_{1}][G_{2}*H][G_{3}]$, the number 
of the branches of $T_{G_{2}}$ is $\lambda=\kappa_{+}=\omega_{3}$ by
Lemma 4.

Let $T$ be any $\omega_{1}$--tree in $M[G_{1}][G_{2}*H][G_{3}]$.
Then there exists $K\subseteq \kappa^{++}$ with $|K|=\omega_{1}$ such that
$T\in M[G_{1}][G_{2}*H][G'_{3}]$, where $G'_{3}=G_{3}\bigcap 
Fn(K,2,\omega_{1})$. 

If $|{\cal B}(T)|=\omega_{4}$ in $M[G_{1}][G_{2}*H][G_{3}]$, then forcing with
$Fn(\kappa^{++}\smallsetminus K,2,\omega_{1})$ will add new branches to
$T$. This implies $T$ is not a Kurepa tree by Lemma 4.

If $|{\cal B}(T)|=\omega_{2}$ in $M[G_{1}][G_{2}*H][G_{3}]$, then by Lemma 4,
$T$ is already a Jech--Kunen tree with $\omega_{2}$--many branches 
in $M[G_{1}][G_{2}*H][G'_{3}]$. Without loss of generality we can assume
that $K=\omega_{1}$. So \[M[G_{1}][G_{2}*H][G'_{3}]\models 
\mbox{ ``There exists
a Jech--Kunen tree with }\omega_{2}\mbox{--many branches''}.\]
But \[M[G_{1}][G_{2}*H][G'_{3}]=M[G'_{3}][G_{1}][G_{2}*H]=
\overline{M}[G_{1}][G_{2}*H],\] where $\overline{M}=M[G'_{3}]$. 
By the same proof of Theorem 7, we can also show that there are no
Jeck--Kunen trees in $\overline{M}[G_{1}][G_{2}*H]$,
a contradiction. \quad $\Box$

\bigskip

Institute of Mathematics, 

The Hebrew University, 

Jerusalem, Israel.

\bigskip

Department of Mathematics, 

Rutgers University, 

New Brunswick, NJ, 08903, USA.

\bigskip

Department of Mathematics,

University of Wisconsin, 

Madison, WI 53706, USA.

\bigskip

{\em Sorting:} The first two addresses are the first author's; the last
one is the second author's. 

\end{document}